\documentstyle{amsppt}
\baselineskip18pt
\magnification=\magstep1
%\NoPageNumbers
%\NoRunningHeads
%\pagewidth{4.5in}
%\pageheight{7.0in}
\pagewidth{30pc}
\pageheight{45pc}
\hyphenation{co-deter-min-ant co-deter-min-ants pa-ra-met-rised
pre-print
fel-low-ship}
\def\leaderfill{\leaders\hbox to 1em{\hss.\hss}\hfill}

\def\H{{\Cal H}}

\def\idest{i.e.\ }
\def\a{{\alpha}}
\def\b{{\beta}}

\def\e{{\varepsilon}}

\def\l{{\lambda}}

\def\t{{\tau}}
\def\w{{\omega}}

\def\extln{{D_n}}
\def\tln{{TL(\widehat{A}_{n-1})}}
\def\otln{{O_n}}
\def\annn{\text{\rm Ann({\bf n})}}
\def\boxit#1{\vbox{\hrule\hbox{\vrule \kern3pt
\vbox{\kern3pt\hbox{#1}\kern3pt}\kern3pt\vrule}\hrule}}
\def\rabbit{\vbox{\hbox{\kern0pt
\vbox{\kern0pt{\hbox{---}}\kern3.5pt}}}}

\def\tableau#1{
	\hbox {
		\hskip -10pt plus0pt minus0pt
		\raise\baselineskip\hbox{
		\offinterlineskip
		\hbox{#1}}
		\hskip0.25em
	}
}

\def\tabCol#1{
\hbox{\vtop{\hrule
\halign{\strut\vrule\hskip0.5em##\hskip0.5em\hfill\vrule\cr\lower0pt
\hbox\bgroup$#1$\egroup \cr}
\hrule
} } \hskip -10.5pt plus0pt minus0pt}

\def\CR{
	$\egroup\cr
	\noalign{\hrule}
	\lower0pt\hbox\bgroup$
}

% Set up the map arrows for commutative diagrams.
%\def\mapright#1{\smash{
%     \mathop{\longrightarrow}\limits^{#1}}}

%Set up macro for commutative diagrams etc. (see Ex. 18.46 in TeXbook)

\def\xfga{{\bf 1}}
\def\xgl{{\bf 2}}
\def\xjona{{\bf 3}}
\def\xms{{\bf 4}}
\def\xmsa{{\bf 5}}
\def\xtl{{\bf 6}}
\def\xwes{{\bf 7}}

\topmatter
\title On Representations of Affine Temperley--Lieb Algebras
\endtitle

\author R.M. Green\endauthor
\affil 
Mathematical Institute\\ Oxford University\\ 24--29 St. Giles'\\
Oxford OX1 3LB\\ England\\ 
{\it  E-mail:} greenr\@maths.ox.ac.uk (now rmg\@colorado.edu)
\endaffil

\abstract
We study the finite-dimensional simple modules, over an algebraically
closed field, of the affine Temperley--Lieb algebra corresponding to the
affine Weyl group of type $A$.  These turn out to be closely related
to the simple modules for a certain $q$-analogue of the annular
algebra of V.F.R. Jones.
\endabstract

\thanks
The author was supported in part by an E.P.S.R.C. postdoctoral
research assistantship.
\endthanks
\endtopmatter

\centerline{\bf This paper was published in 1998 in the Proceedings of the}
\centerline{\bf 8th International Conference on Representations of Algebras}

\head 1. Introduction \endhead

The Temperley--Lieb algebra is a finite-dimensional algebra which was 
introduced in [\xtl] and has been extensively studied in many papers,
for example [\xwes].

In [\xfga], the author and C.K. Fan introduced a diagram calculus for
the ``affine'' Temperley--Lieb algebra, an infinite-dimensional
algebra which is related to affine Weyl group $W(\widehat A_l)$ of
type $A$ in the same way as the (ordinary) Temperley--Lieb algebras
are related to the group algebras of the symmetric groups.  One would
like to be able to classify the finite-dimensional irreducible
representations for this algebra, partly because this gives
irreducible representations of the affine Hecke algebra $\H(\widehat
A_l)$ and of the corresponding Weyl group $W(\widehat A_l)$.  The
``even rank'' case (where $l+1$, the number of nodes in the extended
Dynkin diagram, is even) is already well understood (see [\xms], 
[\xmsa]).  

The diagram calculus in [\xfga] is reminiscent of the diagram calculus
of certain ``annular'' algebras which were introduced in [\xjona], and
which (following [\xgl]) we call Jones algebras.  It is natural to
wonder whether there is a strong connection between them, especially
as both are essentially quotients of the affine Hecke algebra.  

This paper sheds light on these problems and explains how they are related.
We will always work in the category of finite-dimensional modules for
the affine Temperley--Lieb algebra, so the term ``irreducible module''
should be understood to mean ``finite-dimensional irreducible
module''.

After reviewing, in \S2, the diagram calculus of [\xfga], we introduce
$q$-analogues of the Jones algebras in \S3.  (The algebras of [\xjona]
give rise to the case $q = 1$.)  These $q$-Jones algebras (and
certain subalgebras of them) turn out to
be cellular algebras in the sense of [\xgl], in the same way as the
ordinary Jones algebras were cellular (see [\xgl, \S6]).  This means
that we can easily classify their irreducible modules over an
algebraically closed field.  In \S4 we explain the connection between
the objects of \S2 and those of \S3.  This means that ``most'' of the
irreducible modules for the affine Temperley--Lieb algebras are the
same as those described in \S3, under certain natural
identifications.

This method deals with all the modules except those in the most
dominant ``cell'' of the algebra in even rank.  These do not behave in
the same way, and we do not deal with them here.  The reader is
referred to [\xms, \S4.3.1], where these modules are considered.

\vskip 20pt

\head 2. Affine Temperley--Lieb Algebras \endhead

The aim of \S2 is to review the definition of the affine
Temperley--Lieb algebra as given in [\xfga, \S4], and to explain how it is
related to certain other infinite-dimensional algebras given by a
calculus of diagrams.

\vskip 20pt
%\eject
\subhead 2.1 The algebra $\extln$ \endsubhead

It turns out for our purposes to be convenient (for various reasons
which will become clear) to work with an algebra of diagrams which
contains the Temperley--Lieb algebra as a proper subalgebra.  We
start this section by defining this bigger algebra.

First, we recall the definition of an affine $n$-diagram from [\xfga,
\S4.1].

\proclaim{Definition 2.1.1}

An affine $n$-diagram, where $n \in {\Bbb Z}$ satisfies $n \geq 3$,
consists of two infinite horizontal rows of nodes lying at the points
$\{ {\Bbb Z} \times \{0, 1\}\}$ of ${\Bbb R} \times {\Bbb R}$,
together with certain curves, called edges, which satisfy 
the following conditions:

\item{\rm (i)}
{Every node is the endpoint of exactly one edge.}
\item{\rm (ii)}
{Any edge lies within the strip ${\Bbb R} \times [0, 1]$.}
\item{\rm (iii)}
{If an edge does not link two nodes then it is an infinite horizontal
line which does not meet any node.  Only finitely many edges are of
this type.}
\item{\rm (iv)}
{No two edges intersect each other.}
\item{\rm (v)}
{An affine $n$-diagram must be invariant under shifting to the left
or to the right by $n$.}
\endproclaim

We will identify any two diagrams which are isotopic to each other, so
that we are only interested in the equivalence classes of affine
$n$-diagrams up to isotopy.  This has the effect that the only
information carried by edges which link two nodes is the pair of
vertices given by the endpoints of the edge.  Later on in the paper we
will consider ``annular involutions'' which allow the diagrams to be
described without the use of curves, but the geometric viewpoint is
often helpful.

Another way to describe the equivalence of diagrams is as follows.
If $D_1$ and $D_2$ are two affine $n$-diagrams and
there is a bijection $\phi$ between the edges of $D_1$ with the edges
of $D_2$ such that the set of endpoints of each edge $E \in D_1$
(possibly empty) is equal to the set of endpoints of the edge
$\phi(E)$ in $D_2$ (possibly empty), 
then we say that $D_1$ and $D_2$ give the same graph,
and we identify the diagrams.

Because of the condition {\rm (v)},
one can also think of affine $n$-diagrams as
diagrams on the surface of a cylinder, or within an annulus, in a natural way.
Unless otherwise specified, we
shall henceforth regard the diagrams as diagrams on the surface of
a cylinder with $n$ nodes on top and $n$ nodes on the bottom.
Later, we will compare this with Jones' (finite-dimensional) annular 
algebra, which was introduced in [\xjona].  From now on,
we will call the affine $n$-diagrams ``diagrams'' for short, when the
context is clear.  Under this construction, the top row of nodes
becomes a circle of $n$ nodes on one face of the cylinder, which we
will refer to as the top circle.  Similarly, the bottom circle of the
cylinder is the image of the bottom row of nodes.

The cylindrical viewpoint allows one to see why there is an infinite
number of affine $n$-diagrams:
edges which connect nodes on the top circle to
nodes on the bottom circle may wind round the cylinder an arbitrary
number of times, and different winding numbers give nonequivalent diagrams.
We note that the idea of diagrams about a cylinder is familiar from
Potts models in statistical mechanics (see the ``boundary diagrams''
in [\xms]).

An example of an affine $n$-diagram for $n = 4$ is given in Figure 1.
The dotted lines denote the periodicity, and should be identified to
regard the diagram as inscribed on a cylinder.

\topcaption{Figure 1} An affine 4-diagram\endcaption
\centerline{
\hbox to 3.638in{
\vbox to 0.888in{\vfill
	\includegraphics{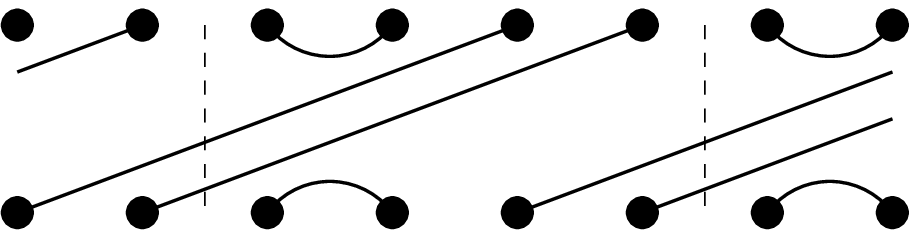}
}
\hfill}
}

\proclaim{Definition 2.1.2}
An edge of the diagram $D$ is said to be vertical if it connects a
point in the top circle of the cylinder to a point in the bottom circle, and
horizontal if it connects two points in the same circle of the cylinder.
\endproclaim

Two diagrams, $A$ and $B$
``multiply'' in the following way, which was described in [\xfga, \S4.2].
Put the cylinder for $A$ on top of the cylinder for $B$ and identify
all the points in the middle row.  This produces a certain (natural) number $x$
of loops.  Removal of these loops forms another diagram $C$ satisfying
the conditions in Definition 2.1.1.  The product $AB$ is then defined
to be $[2]^x C$, where $[2]$ is the Laurent polynomial $v + v^{-1}$ in
the indeterminate $v$. 
It is clear that this defines an associative multiplication.

\proclaim{Definition 2.1.3}
Let $R$ be a ring containing an invertible
indeterminate $v$.  We define the associative 
algebra $\extln$ over $R$ to be the $R$-linear span of all the affine 
$n$-diagrams, with multiplication given as above.
\endproclaim

\vskip 20pt
%\eject
\subhead 2.2 The affine Temperley--Lieb algebra \endsubhead

We assume from now on that $n \geq 3$, so there are at least three
nodes in each circle of the cylinder.

\proclaim{Definition 2.2.1}
Let $R$ be a commutative ring (with identity) containing $[2]$.  
The affine Temperley--Lieb algebra $\tln$ is the $R$-algebra 
given by generators $E_1, \ldots, E_n$ and defining relations $$\eqalignno{
E_i^2 &= [2] E_i, & (1)\cr
E_i E_j &= E_j E_i, \quad \text{ if $\bar{i} \ne \overline{j \pm
1}$}, & (2) \cr
E_i E_{\overline{i\pm 1}} E_i &= E_i. & (3)\cr
}$$  Here, $\bar{i}$ denotes the congruence class modulo $n$ of $i$.
\endproclaim

\proclaim{Remark 2.2.2}\rm
If $n$ is even, and $R$ is a field of characteristic $0$,
this agrees with the ``periodic Temperley--Lieb
algebra'' studied in [\xms], although our notation is different.  The
generator $U_i$ in [\xms] corresponds to our generator $E_{2i-1}$, and
the generator $U_{i j}$ (where $j = i+1$) 
in [\xms] corresponds to our generator $E_{2i}$.

The algebra $\tln$ occurs as a quotient of the affine Hecke algebra of
type $\widehat A$ with algebra generators $T_1, T_2, \ldots T_n$.  The
generator $E_i$ is the image of $v^{-1} T_e + v^{-1} T_i$, where $T_e$
is the identity in the Hecke algebra, and the parameter $q$ in the
Hecke algebra is identified with $v^2$.
\endproclaim

It was shown in [\xfga, \S4] that $\tln$ can also be realised as an
algebra of diagrams, as follows.

\proclaim{Proposition 2.2.3}
The algebra $\tln$ is the subalgebra of $\extln$ spanned by diagrams, $D$,
with the following additional properties: 

\item{\rm (i)}
{If $D$ has no horizontal edges, then $D$ is the identity diagram, in
which point $j$ in the top circle of the cylinder is connected to point
$j$ in the bottom circle for all $j$.}
\item{\rm (ii)}
{If $D$ has at least one horizontal edge, then the number of
intersections of $D$ with the line $x = i + 1/2$ for any integer $i$
is an even number.}

The identification sends the generator $E_i$ to the diagram where
there is a horizontal edge of minimal length connecting $i$ and 
$\overline{i+1}$ in each of the circles of the cylinder, and there is a
vertical edge connecting $j$ in the top circle to $j$ in the bottom circle
whenever $j \ne i, \overline{i+1}$.
\endproclaim

\demo{Proof}
This follows from [\xfga, Corollary 4.5.1].
\qed\enddemo

\demo{Remarks}
The number of intersections of $D$ with a line $x = i + 1/2$ is,
strictly speaking, not well defined, because edges may ``loop back on
themselves''.  We define the ``number of intersections'' as (in (ii))
to be the minimal number possible, which is obtained by perturbing the
diagram of $D$ so that each edge crosses the line as few times as
possible.  Precise details may be found in [\xfga, \S4.4].

A typical example of the representation of a generator $E_i$ as a
diagram on a cylinder is shown in Figure 2.
\enddemo

\topcaption{Figure 2} The element $E_2$ for $n = 5$ \endcaption
\centerline{
\hbox to 1.652in{
\vbox to 1.805in{\vfill
	\includegraphics{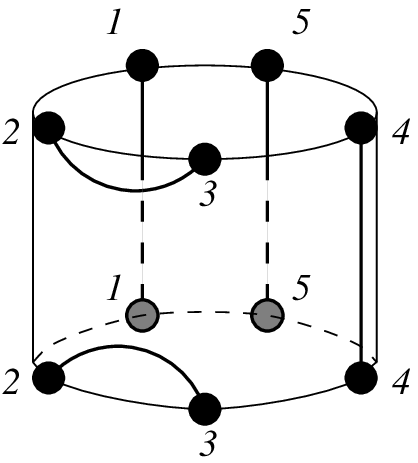}
}
\hfill}
}

We now define an important element $u \in \extln$, similar to the
element $u$ of Jones' annular algebra which was considered in [\xjona].

\proclaim{Definition 2.2.4}
The diagram $u$ of $\extln$ is the one satisfying the property that for
all $j \in \text{\bf n}$, the point $j$ in the bottom circle is connected
to point $\overline{j+1}$ in the top circle by a vertical edge taking
the shortest possible route.
\endproclaim

In the case $n = 4$, the element $u$ is as shown in Figure 3.

\topcaption{Figure 3} The element $u$ for $n = 4$ \endcaption
\centerline{
\hbox to 3.638in{
\vbox to 0.888in{\vfill
	\includegraphics{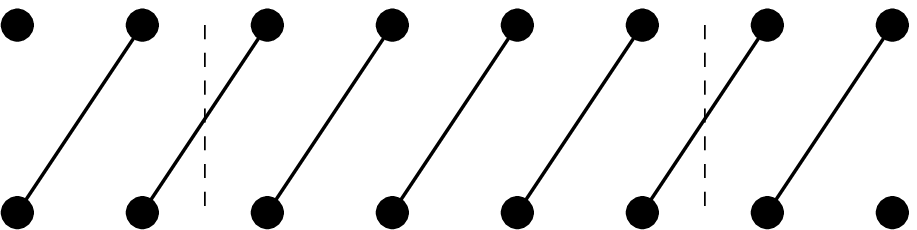}
}
\hfill}
}

\vskip 20pt

\subhead 2.3 The algebra $\otln$ \endsubhead

We also introduce a subalgebra $\otln$ of $\extln$.  This is related
to $\extln$ in the same way as the alternating groups are related to
the symmetric groups.  (The letter ``O'' conceptually stands for
``oriented'', which is terminology borrowed from the paper [\xjona].)

Recall from the remarks following Proposition 2.2.3 that the number of
intersections of a diagram $D$ with a line $x = i + 1/2$ refers to the
minimum possible such number (since we are allowed to replace $D$ by one
of its isotopic images).  In the results which follow, we are usually
only interested in whether this number is even or odd, in which case
it is not necessary that we take the minimum possible number of
intersections, provided that any edges tangent to the line $x = i +
1/2$ are not counted as intersections.

\proclaim{Definition 2.3.1}
Define $\otln$ to be the $R$-submodule of $\extln$ spanned by diagrams $D$
such that the number of
intersections of $D$ with the line $x = i + 1/2$ for any integer $i$
is an even number.

Such diagrams are said to have the even intersection property.
Conversely, if all the numbers of intersections of a diagram $D$ with
the lines $x + 1/2$ are odd, $D$ is said to have the odd intersection property.
\endproclaim

\proclaim{Lemma 2.3.2}
The following identity holds in $\extln$, regarding the generators
$E_i$ of $\tln$ as elements of $\extln$: $$
\underbrace{
E_{\overline{i+2}} E_{\overline{i+3}} \cdots 
E_{\overline{i}}}_{n-1 \text{\ terms}} = u^2 E_i
.$$
\endproclaim

\demo{Proof}
Since $u E_i u^{-1} = E_{\overline{i+1}}$, it suffices to check the
case $i = 1$.  This follows easily from the diagrammatic viewpoint.
\qed\enddemo

\proclaim{Lemma 2.3.3}
The $R$-module $\otln$ is the subalgebra of $\extln$ generated by
$\tln$ and the elements $\{u^{2m} : m \in {\Bbb Z}\}$.
\endproclaim

\demo{Proof}
It is not hard to see that if a diagram has no horizontal edges,
it must be a power of $u$, and that the diagram lies in
$\otln$ if and only if this power of $u$ is even.  
The fact that there are no more diagrams in $\otln$ now follows from
Proposition 2.2.3.

To check that these diagrams span a subalgebra, it remains to
prove that $u^2 E_i$ and $u^{-2} E_i$ lie in $\otln$
for each $i$.  The case of $u^2 E_i$ follows from Lemma 2.3.2, and
the case of $u^{-2} E_i$ follows by an analogous argument.
The assertion now follows.
\qed\enddemo

The following lemma shows that any diagram has the even or odd
intersection property, according as it is an element of $\otln$ or
not.

\proclaim{Proposition 2.3.4}
Let $D$ be a diagram in $\extln$.  If $D \in \otln$ (respectively, $D
\not\in \otln)$, then the number
of intersections of $D$ with the line $x = i + 1/2$ is an even
(respectively, odd) number for each $i$.
\endproclaim

\demo{Proof}
If $D \in \otln$, it is clear from Lemma 2.3.3 that $D$ has the
even intersection property, 
because it is
true if $D$ is an even power of $u$ or an element of $\tln$.  The
converse, that any diagram with the even intersection property lies in
$\otln$, is equally clear.

Next, we observe that if $D$ has the even (respectively, odd)
intersection property, then $u.D$ has the odd (repsectively, even)
intersection property.  This follows by consideration of the
product $u.D$ via diagrams.

By consideration of the fact that each point in the diagram $D$ has exactly
one edge emerging from it, and that this edge can either emerge
towards the left, towards the right, or straight upwards, we see that
any diagram $D$ in $\extln$ has either the odd intersection property
or the even intersection property as mentioned in the statement.

The proof now follows.
\qed\enddemo

\proclaim{Lemma 2.3.5}
As $O_n$-bimodules, $\extln = \otln \oplus u \otln = \otln \oplus \otln u$.
\endproclaim

\demo{Proof}
We observe from the proof of Proposition 2.3.4 that
any diagram of $\extln$ is either a diagram of $\otln$ or of the form
$u . D$ where $D$ is a diagram of $\otln$.  We also note that 
for a diagram $D$, $D \in \otln \Rightarrow u^{-1} D u \in \otln$.
The proof follows from these facts.
\qed\enddemo

We now define an automorphism $\e$ of $\extln$ which fixes $\otln$.
%\vfill\eject
\proclaim{Proposition 2.3.6}
The $R$-linear map $\e : \extln \rightarrow \extln$ is defined by its
effect on the diagrams $D$ as follows: $$\e(D) = \cases
D & \text{ if } D \in \otln;\cr
-D & \text{ if } D \not\in \otln.\cr
\endcases$$  This map is an automorphism of $\extln$ whose fixed points are
the elements of $\otln$.

Thus we have $\e(u) = -u$, $\e(u^{-1}) = -u^{-1}$ and $\e(E_i) =
 E_i$ for all $i$.
\endproclaim

\demo{Proof}
From Lemma 2.3.5, we find that any diagram in $\extln
\backslash \otln$ can be written (uniquely) in the form $u . E$
where $E$ is a diagram of $\otln$.  
Conversely, elements of these forms do not lie
in $\otln$.  Using the fact that $u . \otln . u^{-1} = \otln$ and the fact
that even powers of $u$ lie in $\otln$, we
see that the map $\e$ as described respects the multiplication
(because two diagrams in $\extln \backslash \otln$ multiply to give an
element of $\otln$, and so on).

The effect of $\e$ on the generators and the assertion about the
fixed points of $\e$ are easy to prove.
\qed\enddemo

To conclude this section, we note that $\extln$ has a nice presentation in
terms of generators and relations, as follows.

\proclaim{Proposition 2.3.7}
The algebra $\extln$ is generated by elements $$E_1, \ldots, E_n, u,
u^{-1}.$$  It is subject to the relations {\rm (1)}, {\rm (2)} and
{\rm (3)} of
Definition 2.2.1, and the following additional defining 
relations: $$\eqalignno{
u E_i u^{-1} &= E_{\overline{i+1}}, & (4)\cr
(u E_1)^{n-1} &= u^n . (u E_1). & (5)\cr
}$$
\endproclaim

\demo{Proof}
From Lemma 2.3.3 and Lemma 2.3.5, we see that the given
set does generate the algebra $\extln$.  In addition, all the
relations are true; (5) looks unfamiliar but we claim that
it is equivalent to Lemma 2.3.2 after applications of (4).  
To do this, first reexpress (5) as $u^{-(n-1)} (u E_1)^{n-1} =
u^2 E_1$.  Now apply (4) to quasi-commute all the occurrences $u$ in
the left hand side to the left of the expression.  This establishes the
special case of Lemma 2.3.2 where $i = 1$.  To obtain the general
expression, we conjugate each side of the equation by a suitable power
of $u$.

We then observe that any word in the generators can be rewritten as a
multiple of a power of $u$ times a monomial in the $E_i$ (\idest a
diagram of $\tln$) by using the
relations (1) to (4).  Using Lemma 2.3.2 (\idest relation (5)), one
can reduce this power of $u$ to $0$ or $1$, unless the word is simply
a power of $u$.  
Lemma 2.3.5 now shows that the relations are sufficient.
\qed\enddemo

\proclaim{Remark 2.3.8}\rm
In fact, $\extln$ is generated by $u$, $u^{-1}$ and $E_1$, as can be
seen from relation (4).
\endproclaim

\vskip 20pt

\head 3. $q$-Jones Algebras \endhead

In \S3, we assume that $K = F(v)$ is an algebraically closed field,
where $F$ is a subfield of $K$ and $v$ is an indeterminate.  We
then introduce a family of finite-dimensional algebras $J_q(n)$, where
$q$ is a nonzero element of $K$.  We call $J_q(n)$ the $q$-Jones
algebra of rank $n$, where $n$ is a natural number at least $3$.  
The behaviour of these algebras in odd and in even ranks turns out to
be distinctly different; if $n$ is odd then $J_1(n)$ is Jones's
annular algebra, $J(n)$, as defined in [\xjona].  However, if $n$ is even,
$J_1(n)$ is not the full Jones algebra but a quotient of the
subalgebra spanned by the ``oriented elements'' of $J(n)$.  In both
even and odd case, the algebras turn out to be cellular in the sense
of [\xgl], which enables their simple modules over $K$ to be
easily determined.  These $q$-Jones algebras appear to be new.

\vskip 10pt

\subhead 3.1 Cellular algebras \endsubhead

For completeness, we recall the definition of a cellular algebra from [\xgl].

\proclaim{Definition 3.1.1}  Let $R$ be a commutative
ring with identity.  A {\it cellular algebra} over $R$ is an associative unital
algebra, $A$, together with a cell datum $(\Lambda, M, C, *)$ where

\item {\rm 1.}
{$\Lambda$ is a poset.  For each $\l \in \Lambda$, $M(\l)$ is a finite set
(the set of ``tableaux'' of type $\l$) such that $$
C : \coprod_{\l \in \Lambda} \left( M(\l) \times M(\l) \right) \rightarrow A
$$ is injective with image an $R$-basis of $A$.}
\item {\rm 2.}
{If $\l \in \Lambda$ and $S, T \in M(\l)$, we write $C(S, T) = C_{S, T}^{\l}
\in A$.  Then $*$ is an $R$-linear involutary anti-automorphism 
of $A$ such that
$(C_{S, T}^{\l})^* = C_{T, S}^{\l}$.}
\item {\rm 3.}
{If $\l \in \Lambda$ and $S, T \in M(\l)$ then for all $a \in A$ we have $$
a . C_{S, T}^{\l} \equiv \sum_{S' \in M(\l)} r_a (S', S) C_{S', T}^{\l}
\mod A(< \l),
$$ where  $r_a (S', S) \in R$ is independent of $T$ and $A(< \l)$ is the
$R$-submodule of $A$ generated by the set $$
\{ C_{S'', T''}^{\mu} : \mu < \l, S'' \in M(\mu), T'' \in M(\mu) \}
.$$}
\endproclaim

\subhead 3.2 Annular Involutions \endsubhead

Let $D$ be a diagram (as in \S2) representing a basis element of the
extended affine Temperley--Lieb algebra $\extln$.  Throughout \S3.2,
we are only concerned with diagrams $D$ with at least one vertical
edge.  We define a certain integer $w(D)$ as follows.

\proclaim{Definition 3.2.1}
Let $D$ be as above.  
Let $w_1(D)$ be the number of pairs $(i, j) \in {\Bbb Z} \times {\Bbb Z}$
where $i > j$ and $\overline{j}$ in the bottom circle of
$D$ is joined to $\overline{i}$ in the top circle of $D$ by an edge
which crosses
the ``seam'' $x = 1/2$.  We then define $w_2(D)$ similarly but
with the condition that $i < j$, and we define $w(D) = w_1(D) - w_2(D)$.
\endproclaim

\proclaim{Remark 3.2.2}\rm
We think of $w(D)$ as measuring some kind of ``winding
number'' of the diagram $D$.  

The non-intersection condition guarantees that at least one of
$w_1(D)$ or $w_2(D)$ is $0$, and that
$w(D)$ is always finite.  However, the values of $w(D)$ for $D
\in \extln$ are unbounded even though $n$ is fixed.

Some examples of the values of $w(D)$ are as follows.  In Figure 1,
the pairs which cross the seam are $(3, -1)$ and $(4, 0)$, giving
$w(D) = 2$.  Figure 2 has $w(D) = 0$.  Figure 3 has $w(D) = 1$,
corresponding to the pair $(1, 0)$.  A moment's thought shows that
$w(u^n) = n$ for any $n \in {\Bbb Z}$, which motivates the definition.
\endproclaim

We recall the definition of an annular involution of the symmetric group
from [\xgl, Lemma 6.2].

\proclaim{Definition 3.2.3}
An involution $S \in {\Cal S}_n$ is annular if and only if, for each
pair $i, j$ interchanged by $S$ $(i < j)$, we have
\item{\rm(a)}
{$S[i, j] = [i, j]$ and}
\item{\rm(b)}
{$[i, j] \cap \text{\rm Fix} S = \emptyset$ or $\text{\rm Fix}S
\subseteq [i, j]$.}

We write $S \in I(t)$ if $S$ has $t$ fixed points, and we write $S \in
\annn$ if $S$ is annular.
\endproclaim

\proclaim{Lemma 3.2.4}
The diagrams of $\extln$ with vertical edges are in canonical
bijection with triples 
$[S_1, S_2, w]$ where $S_1, S_2 \in \annn \cap I(t)$ for $t > 0$, and 
$w \in {\Bbb Z}$.  
\endproclaim

\demo{Proof}
The correspondence is given as follows.  Consider a diagram $D$.
The involution $S_1$
exchanges points $i$ and $j$ if and only if $i$ is connected to $j$ in
the top circle of $D$.  Similarly $S_2$ exchanges points $i$ and $j$ if
and only if $i$ is connected to $j$ in the bottom circle of $D$.  One
can easily see that $S_1$ and $S_2$ are annular.
Finally $w = w(D)$; this determines which points are connected by the
vertical edges (it is enough to specify one vertical edge, by the
non-intersection criterion).

Conversely consider a triple $[S_1, S_2, w]$ as required.  The
hypothesis $t > 0$ guarantees
that each of $S_1$ and $S_2$ each has a fixed point.  If $S \in
{S_1, S_2}$ exchanges points $i$ and $j$, it now follows from  part (b)
of the definition that we can determine which way round the cylinder
the corresponding edge in the diagram goes.  Using part (a) we can
show that no intersections occur, so $S_1$ determines the horizontal
edges in the top of the diagram and $S_2$ determines those in the
bottom.  As regards $w$, observe that one can disconnect all the
vertical edges from a diagram and reconnect the top end of
each one to the next available point to the right.  This increases
$w(D)$ by 1.  Conversely one can decrease it by one.  Thus $w$
determines the endpoints of all the vertical edges, which completes the proof.
\qed\enddemo

\proclaim{Definition 3.2.5}
For $S_1, S_2 \in \annn \cap I(t)$ where $t > 0$, 
we define $r(S_1, S_2)$ to be the smallest nonnegative integer (which
will be $0$ or $1$, following Proposition 2.3.4) satisfying $[S_1,
S_2, r(S_1, S_2)] \in \otln$.
\endproclaim

\vskip 20pt
%\vfill\eject
\subhead 3.3 $q$-Jones Algebras in Odd Rank \endsubhead

In \S3.3, we assume that $n$ is an odd number.  (We deal with
this before the even case because the even case is already better understood.)

We now put an equivalence relation, $\sim$, on the basis elements 
$[S_1, S_2, w]$.

\proclaim{Definition 3.3.1}
We say that $[S_1, S_2, w] \sim [T_1, T_2, w']$ if and only if $S_1 =
T_1$, $S_2 = T_2$ and $w$ and $w'$ are congruent modulo $t$, where
$S_1, S_2 \in \annn \cap I(t)$.  

We say a basis element $[S_1, S_2, w]$ is in root position if and only if 
$0 \leq w < t$.
Clearly the set of basis elements in root position is a set of
equivalence class representatives for $\sim$.  We usually
denote such an element by $[S_1, S_2, \t_i]$ to emphasize that it is
in root position.
\endproclaim

\proclaim{Lemma 3.3.2}
The subspace spanned by elements of the form $$
\{ [S_1, S_2, w] - q^s [T_1, T_2, w'] : [S_1, S_2, w] \sim [T_1, T_2,
w'] \text{\rm \ and \  } w' = w + ts \},
$$ where $S_1 \in \annn \cap I(t)$, is an ideal of $\extln$.  We 
denote it by $\w(q)$.
\endproclaim

\demo{Proof}
It is clear from the diagram calculus that the
product of any two basis elements of $\extln$ is a multiple of another one.
If $$
[X_1, X_2, w] [Y_1, Y_2, w'] = \b [Z_1, Z_2, w'']
$$ for some scalar $\b$ then one can also easily see that $$
[X_1, X_2, w] [Y_1, Y_2, w' + t] = \b [Z_1, Z_2, w'' + t].$$
The claim now follows.
\qed\enddemo

\proclaim{Definition 3.3.3}
The $q$-Jones algebra over $K(q)$ (where $q$ is an indeterminate) is
the quotient of $\extln$ by $\w(q)$.
\endproclaim

\proclaim{Remark 3.3.4}\rm
When $q = 1$, one recovers the usual Jones algebra (see [\xgl,
Proposition 6.14]).
\endproclaim

\proclaim{Lemma 3.3.5}
The $q$-Jones algebra has a basis given by the images of the elements $$
\{ [S_1, S_2, \t_i] : S_1, S_2 \in \annn \cap I(t), \t_i \in
\{0, 1, \ldots, t-1\} \}
$$
\endproclaim

\demo{Proof}
It is clear that the given elements span the algebra.

Suppose for a contradiction that there is some nontrivial linear
relation between the elements.  Then we may assume that the
coefficients of this relation lie in $K[q]$ and that not all of them have
$(q-1)$ as a factor.  Putting $q = 1$ and using Remark 3.3.4, we now
obtain a nontrivial linear relation for the corresponding basis
elements of the Jones algebra, which is a contradiction.
\qed\enddemo

\vskip 20pt

\subhead 3.4 
Representation theory of the $q$-Jones algebra in odd rank \endsubhead

We can now prove that the $q$-Jones algebra is cellular.  The
motivation behind doing this is that we will then have bases for all
its finite dimensional irreducible modules over $K$.  The methods are the same
as those in [\xgl, \S6].  

We specialize the
indeterminate $q$ to a fixed, nonzero element of $K$.  We also
define the set ${\Cal T}(n)$ to be the set of odd natural
numbers up to and including $n$.  We abuse notation by identifying an
element $[S_1, S_2, w]$ where $0 \leq w < t$ with its image in $J_q(n)$.

\proclaim{Definition 3.4.1}
Let $w_0$ be the longest element of ${\Cal S}_n$.  If $S \in \annn$,
we define $S^* \in \annn$ to be the element $w_0^{-1} S w_0$.
\endproclaim

\proclaim{Proposition 3.4.2}
Let $t \in {\Cal T}(n)$.  We factorise $$
x^t - q = \prod_{i = 1}^t (x - r_i(t, q))
$$ and write $$
f_j(x, q) = \prod_{i > j}(x - r_i(t, q)) = \sum_{i = 0}^{t - j}r_{ij}(t, q) x^i
$$ for $j = 1, \ldots, t$.  Write $f_0(x, q) = 1$ so that $r_{00}(0) =
1$.  Then a cell datum for $J_q(n)$ over $K$ is given by $(\Lambda, M,
C, *)$ defined as follows: 

\item{\rm(a)}
{$\Lambda = \{(t, j): t \in {\Cal T}(n), j
\in \text{\bf t}\}$, ordered lexicographically, \idest $(t', j') < (t,
j)$ if and only if $t' < t$ or $t = t', j' < j$.}

\item{\rm(b)}
{For $(t, j) \in \Lambda$, $M(t, j) = \annn \cap I(t)$.}

\item{\rm(c)}
{If $S_1, S_2 \in M(t, j)$ then $C_{S_1, S_2}^{t, j} = 
\sum_{i = 0}^{t - j} r_{ij}(t, q) [S_1, S_2^*, i].$}

\item{\rm(d)}
{The anti-automorphism, $*$, is defined by $[S_1, S_2, w]^* = 
[S_2^*, S_1^*, w]$.}
\endproclaim

\demo{Proof}
We check that the three conditions of Definition 3.1.1 hold.
Condition 1 comes from Lemma 3.3.5.

To prove condition 2, we observe that applying $*$ to a
basis element $[S_1, S_2^*, w]$ in root position gives another element
in root position.  The fact that the map is an anti-automorphism is
clear from the diagrammatic viewpoint---the cylinders are being
rotated so they are upside-down.

The proof of condition 3 is the same as [\xgl, Theorem 6.15], with a few
minor and obvious changes.
\qed\enddemo

The motivation behind proving Proposition 3.4.2 is that the theory of
cellular algebras allows one to construct a complete list of 
finite-dimensional irreducible modules for $J_q(n)$.

\proclaim{Definition 3.4.3}
Let $\l = (t, j) \in \Lambda$.
The module $W(t, j, q)$ for $J_q(n)$ has a basis given by $\{C_S :
S \in M(t, j)\}$ with action specified as follows.  If $a \in J_q(n)$
and $$
a . C_{S, T}^{t, j} \equiv \sum_{S' \in M(t, j)} c_a (S', S) C_{S',
T}^{t, j}
\mod J_q(n)\ (< \l)
$$ for any suitable $T$ (for example, $S$), then the left action of
$a$ on $W(t, j, q)$ is defined by $$
a . C_S = \sum_{S' \in M(t, j)} c_a (S', S) C_{S'}
.$$
\endproclaim

\proclaim{Remark 3.4.4}\rm
The theory of cellular algebras in [\xgl] ensures that this is
well-defined.
\endproclaim

\proclaim{Proposition 3.4.5}
The modules $W(t, j, q)$, as $r_j(t, q)$ ranges over the $t$-th roots of $q$,
are a complete set of irreducible representations for $J_q(n)$ over $K$.
\endproclaim

\demo{Proof}
This is essentially a restatement of [\xgl, Corollary 6.17], after making
the necessary trivial changes to [\xgl, Scholium 6.16].  Note that the
parameter $\delta$ occurring in [\xgl, Corollary 6.17] is equal to
$[2] = v + v^{-1} \ne 0$ in our case.
\qed\enddemo

\proclaim{Proposition 3.4.6}
Let $\psi$ be a representation corresponding to the irreducible module $W(t, j,
q)$, and let $[S_1, S_2, w]$ be a basis element of $J_q(n)$
arising from the set $M(t, j)$.  Then $$
\psi([S_1, S_2, w + 1]) = r_j(t, q) . \psi([S_1, S_2, w])
.$$
\endproclaim

\demo{Proof}
This is a corollary of the proof of Proposition 3.4.2.  For more
details see [\xgl, (6.15.1), (6.15.2)].
\qed\enddemo

\proclaim{Remark 3.4.7}\rm
Proposition 3.4.6 gives a convenient way to test whether the
representations $W(t, j, q)$ and $W(t, j', q)$ are distinct.  We will
use this fact later in the classification of the simple modules for $\tln$.
\endproclaim

\vskip 20pt

\subhead 3.5 $q$-Jones Algebras in Even Rank \endsubhead

We now highlight the differences between the odd rank situation and
the corresponding even rank situation.  We treat the case of even rank
in less detail
than the case of odd rank studied in \S3.3 and \S3.4, because the
situation is already better understood and the techniques are very similar.

In \S3.5, $n$ is even and the algebra $\otln$ takes over the role which
$\extln$ played in \S3.3 and \S3.4.  
Since $n$ is even, any diagram has an even number of vertical edges,
and so the relation $\sim$ of
Definition 3.3.1 carries over naturally to an equivalence relation on
the basis elements of $\otln$.  

For $n$ even, we set ${\Cal T}(n) := \{ 0, 2, 4, \ldots, n \}$.

\proclaim{Lemma 3.5.1}
The diagrams of $\otln$ with vertical edges are in canonical bijection
with triples 
$[S_1, S_2, r(S_1, S_2) + w]$ where $t \in {\Cal T}(n)\backslash \{0\}$,
$S_1, S_2 \in \annn \cap I(t)$, and $w \in 2{\Bbb Z}$.
\endproclaim

\demo{Proof}
This is along the same lines as the proof of Lemma 3.2.4.
\qed\enddemo

\proclaim{Remark 3.5.2}\rm
The hypothesis that there exists a vertical edge in the diagram is
necessary, otherwise one could form a counterexample using any diagram
with no vertical edges!
\endproclaim

\proclaim{Definition 3.5.3}
The $R$-module $\w'(q)$ of $\otln$ is the span of elements of $\otln$ of the
form given in Lemma 3.3.2 where $[S_1, S_2, w]$, and hence $[T_1, T_2,
w']$, involve vertical edges.  The ideal
$I_0$ of $\otln$ is the span of the elements of $\otln$ with no vertical edges.

The $q$-Jones algebra over $K(q)$ (where $q$ is an indeterminate) is
the quotient of $\otln$ by the ideal $\w'(q) + I_0(q)$.
\endproclaim

\proclaim{Remark 3.5.4}\rm
It is clear from the diagrammatic viewpoint that $I_0(q)$ is an ideal.
One proves that $\w'(q) + I_0(q)$ is an ideal by using an argument identical to
that of the proof of Lemma 3.3.2.  

\vskip 10pt

It should be noted that the term ``$q$-Jones algebra'' is
misleading in even rank, because specializing $q$ to $1$ gives a
smaller algebra than Jones' annular algebra.  In fact, we obtain the
quotient of the Jones algebra by its ideal $J_0(n)$ spanned by basis elements
with no through-strings.  The reason for factoring out $I_0$ is that 
irreducible modules $M$ satisfying $I_0.M = M$ behave rather differently
from those satisfying $I_0.M = 0$.  In this paper we are interested in
the latter class of modules.
\endproclaim

\proclaim{Lemma 3.5.5}
The $q$-Jones algebra $J_q(n)$, for even $n$,
has a basis given by the images of the elements $$
\{ [S_1, S_2, r(S_1, S_2) + 2w] \}
,$$ where
$t \in {\Cal T}(n)\backslash \{0\}$, $S_1, S_2 \in \annn \cap I(t)$
and $w \in \{0, 1, \ldots, t/2 - 1\}.$
\endproclaim

\demo{Proof}
As in Lemma 3.3.5, it is clear that the given elements span the
algebra.

Suppose for a contradiction that there is some nontrivial linear
relation between the elements.  Then we may assume that the
coefficients of this relation lie in $K[q]$ and that not all of them have
$(q-1)$ as a factor.  Putting $q = 1$, we now
obtain a nontrivial linear relation for the corresponding basis
elements of $J(n)/J_0(n)$, and thus a nontrivial linear relation in
$J(n)$, which is a contradiction.
\qed\enddemo

We will abuse notation again by writing $[S_1, S_2, w]$ (one of
the elements appearing in Lemma 3.5.5) for its image in $J_q(n)$.

We now give a cell datum for $J_q(n)$ in even rank.

\proclaim{Proposition 3.5.6}
Define ${\Cal T}'(n) := \{1, 2, \ldots, n/2\}$, and let $t \in {\Cal
T}'(n)$.  We factorise $$
x^{t} - q = \prod_{i = 1}^t (x - r_i(t, q))
$$ and write $$
f_j(x, q) = \prod_{i > j}(x - r_i(t, q)) = \sum_{i = 0}^{t - j}r_{ij}(t, q) x^i
$$ for $j = 1, \ldots, t$.  Write $f_0(x, q) = 1$ so that $r_{00}(0) =
1$.  Then a cell datum for $J_q(n)$ over $K$ is given by $(\Lambda, M,
C, *)$ defined as follows: 

\item{\rm(a)}
{$\Lambda = \{(t, j): t \in {\Cal T}'(n), j
\in \text{\bf t}\}$, ordered lexicographically, \idest $(t', j') < (t,
j)$ if and only if $t' < t$ or $t = t', j' < j$.}

\item{\rm(b)}
{For $(t, j) \in \Lambda$, $M(t, j) = \annn \cap I(2t)$.}

\item{\rm(c)}
{If $S_1, S_2 \in M(t, j)$ then $C_{S_1, S_2}^{t, j} = 
\sum_{i = 0}^{t - j} r_{ij}(t, q) [S_1, S_2^*, r(S_1, S_2^*) + 2i].$}

\item{\rm(d)}
{The anti-automorphism, $*$, is defined by $[S_1, S_2, w]^* = 
[S_2^*, S_1^*, w]$.}
\endproclaim

\demo{Proof}
This uses exactly the same techniques as Proposition 3.4.2.
\qed\enddemo

\proclaim{Remark 3.5.7}\rm
Using this cell datum, one can state and prove results exactly
analogous to Definition 3.4.3, Remark 3.4.4, Proposition 3.4.5,
Proposition 3.4.6 and Remark 3.4.7.  The only appreciable differences
are that the role of ${\Cal T}(n)$ in the odd rank is played by ${\Cal
T}'(n)$ in the even rank, and that the role of $t$ in the odd rank is
played by $t/2$ in the even rank.
\endproclaim

\vskip 20pt

\head 4. Classification of Simple Modules for $\tln$\endhead

To classify the simple modules for $\tln$, we find
once again that it is helpful to consider the cases of odd and even
$n$ separately.  We tackle the case of odd rank first, because it is
less well understood.  
%Throughout \S4, ``simple'' always
%means ``finite-dimensional simple''.

\subhead 4.1 Simple modules in odd rank \endsubhead

The technique in odd rank is to calculate the simple modules for
$\extln$ and study their behaviour upon restriction to $\otln$ and
then to $\tln$.  We also recall the automorphism $\e$ from \S2.

\proclaim{Proposition 4.1.1}
The irreducible representations $\psi$ 
for $\extln$ are parametrised by pairs $(\a_t, t)$,
where $t \in {\Cal T}(n)$ and $\a_t \in K \backslash 0$.  The scalar
$\a_t$ is that appearing in the equation $$
\psi([S_1, S_2, w + 1]) = \a_t . \psi([S_1, S_2, w])
,$$ where $S_1, S_2 \in \annn \cap I(t)$.
\endproclaim

\demo{Proof}
Let $M$ be a simple module for $\extln$.  Let $I_t$ be the ideal of
$\extln$ spanned by all diagrams with at most $t$ vertical edges.
Then there exists a unique (odd) number 
$t$ such that $I_t . M = M$ but $I_{t-2} . M = 0$, interpreting $I_{-1} = 0$.
We also see that $M$ is an irreducible module for $D_n /
{I_{t-2}}$.  Clearly, $t \in {\Cal T}(n)$.

One sees from the diagrammatic viewpoint that $u^n$ lies in the centre
of $\extln$.  By Schur's Lemma, $u^n$ acts as a scalar, $\a$, on
$M$.  The ideal generated by $u^n - \a$ is the same as the ideal
$\w(q)$ in Definition 3.3.3, for $q = \a$, so $M$ is a simple module
for the $q$-Jones algebra $J_q(n)$ for this value of $q$.  Proposition
3.4.5 now shows that the irreducibles are parametrised by the set of
$t$-th roots of $\a$.  (Note that we allow the possibility that $\a$
does not have $t$ distinct $t$-th roots and there are fewer than $t$
irreducibles.)
The fact that nonisomorphic irreducibles have
different parameters comes from Proposition 3.4.6 and Remark 3.4.7.
Conversely, any nonzero $\a \in K$ turns up as one
of the numbers $r_j(t, q)$ for $q = \a^t$.  

This means that all the
irreducibles for $\extln$ are parametrised by pairs $(\a_t, t)$, where
$t$ is as above and $\a_t \in K \backslash 0$.  The last assertion
follows from Proposition 3.4.6.
\qed\enddemo

\proclaim{Proposition 4.1.2}
The irreducible representations $\rho$ 
for $\otln$ are parametrised by pairs $(\a_t, t)$,
where $t \in {\Cal T}(n)$ and $\a_t \in K \backslash 0$.  The scalar
$\a_t$ is that appearing in the equation $$
\rho([S_1, S_2, w + 2]) = \a_t . \rho([S_1, S_2, w])
,$$ where $S_1, S_2 \in \annn \cap I(t)$.
\endproclaim

\demo{Proof}
Let $L$ be a simple module for $\extln$, and let $L'$ be a simple
$\otln$-submodule of $L$ (using the $\otln$-module structure suggested
by Lemma 2.3.5).  Since $u^n$ acts as a scalar on $L$, it
acts as a scalar on $L'$, and in particular $u^n L' = L'$.  However,
$u^2 L' = L'$ since $u^2 \in \otln$.  Since $n$ is odd and $u$ is
invertible, this implies $u . L' = L'$, and thus $L = L'$.
Furthermore, $L$ is determined by $L'$ and the scalar by which $u^n$
acts.

It is clear that two irreducible modules $L_1$ and $L_2$ for $\extln$
restrict to nonisomorphic modules for $\otln$ unless their
parameters $(\a_1,
t_1)$ and $(\a_2, t_2)$ satisfy $t_1 = t_2$ and $\a_1 = \pm \a_2$.
(This is because $$
\psi([S_1, S_2, w + 2]) = \a^2 . \psi([S_1, S_2, w]),
$$ where $\psi$ is a representation of $\extln$ 
corresponding to a pair $(\a, t)$ and $S_1, S_2 \in \annn \cap I(t)$.)  If
$\text{\rm char}K = 2$ then necessarily $L_1 = L_2$.  Otherwise,
let $\psi_1$ be a representation affording $L_1$.  We claim
that if $L_1$ is given by the pair $(\a_t, t)$, then $\psi_1 . \e$
is an irreducible representation for $\extln$ corresponding to the pair
$(-\a_t, t)$ and thus affording $L_2 \ne L_1$.  This is because $$
\psi([S_1, S_2, w + 1]) = \a . \psi([S_1, S_2, w]) \Rightarrow
\psi(\e([S_1, S_2, w + 1])) = - \a . \psi(\e([S_1, S_2, w]))
.$$  Irreducibility follows from the fact that $L_1$ and $L_2$ have
the same dimension.  
Thus $\psi_1 . \e$ affords $L_2$, although $\psi_1$ and $\psi_1
. \e$ give the same representation of $\otln$, since $\e$ fixes $\otln$.

The result now follows from these observations and Frobenius reciprocity.
\qed\enddemo

\proclaim{Theorem 4.1.3}
The nontrivial irreducible representations $\rho$ 
for $\tln$ are parametrised by pairs $(\a_t, t)$,
where $t \in {\Cal T}(n)\backslash \{n\}$ and $\a_t \in K \backslash
0$.  The scalar $\a_t$ is that appearing in the equation $$
\rho([S_1, S_2, w + 2]) = \a_t . \rho([S_1, S_2, w])
,$$ where $S_1, S_2 \in \annn \cap I(t)$.  The only other irreducible
representation is the ``trivial'' representation 
of dimension $1$ which sends the identity to 
$1$ and all the generators $E_i$ to $0$.
\endproclaim

\demo{Proof}
Recall that $\tln$ and $\otln$ are filtered by two-sided ideals
indexed by the number of vertical edges in the diagrams.  Apart from
the top section, which contains only the identity in $\tln$ but is
infinite in $\otln$, the sections of this filtration by ideals are
the same (in fact, the ideals themselves are the same).  It follows
that the simple modules are the same, except for
those corresponding to the top section.  In the latter case, only the
``trivial'' module (as described in the statement of the Theorem) for
$\tln$ survives.
\qed\enddemo

\vskip 20pt

\subhead 4.2 Simple modules in even rank \endsubhead

Similar techniques can also be used to classify the simple modules for
$\tln$ where $n$ is even.  As usual we exclude the case where $I_0
. M = M$ and $I_0$ is as in Definition 3.5.3.

\proclaim{Proposition 4.2.1}
The irreducible representations $\rho$ 
for $\otln$ satisfying $\rho(I_0) = 0$ are parametrised by pairs $(\a_t, t)$,
where $t \in {\Cal T}'(n)$ and $\a_t \in K \backslash 0$.  The scalar
$\a_t$ is that appearing in the equation $$
\rho([S_1, S_2, w + 2]) = \a_t . \rho([S_1, S_2, w])
,$$ where $S_1, S_2 \in \annn \cap I(2t)$.
\endproclaim

\demo{Proof}
Let $M$ be the simple module corresponding to $\rho$.  Again, $u^n$
acts as a scalar, $\a$, on $M$, by Schur's lemma.  Thus the ideal
generated by $u^n - \a$ factors through $M$, and we obtain an
irreducible representation for $J_q(n)$ with $q = \a$.  The list of
such representations is indexed by the $t$-th roots of $\a$, and they
are distinguishable as in the statement of the proposition by the
usual techniques.  As in Proposition 4.1.1, such a representation
exists for any $\a_t$, which completes the proof.
\qed\enddemo

Finally we rederive the result in [\xms] which parametrises the
simple modules for $\tln$ in even rank.

\proclaim{Theorem 4.2.2}
The nontrivial irreducible representations $\rho$ 
for $\tln$ which satisfy $\rho(I_0) = 0$ are parametrised by pairs $(\a_t, t)$,
where $t \in {\Cal T}'(n)\backslash \{n/2\}$ and $\a_t \in K \backslash
0$.  The scalar $\a_t$ is that appearing in the equation $$
\rho([S_1, S_2, w + 2]) = \a_t . \rho([S_1, S_2, w])
,$$ where $S_1, S_2 \in \annn \cap I(2t)$.  The only other irreducible
representation is the ``trivial'' representation
of dimension $1$ which sends the identity to
$1$ and all the generators $E_i$ to $0$.
\endproclaim

\demo{Proof}
This follows the same lines as Theorem 4.1.3, except that it uses
Proposition 4.2.1 instead of Proposition 4.1.2.
\qed\enddemo

\vskip 20pt

\head Acknowledgements \endhead

The author would like to thank C.K. Fan for some helpful discussions.
The author is also grateful to K. Erdmann and the referee
for numerous helpful comments on the manuscript.

\vskip 20pt

\vskip 1cm

\leftheadtext{}
\rightheadtext{}
%\vfill\eject

\end